\documentclass[preprint,3p,sort]{elsarticle}

\usepackage{amsfonts}
\usepackage{amsthm}
\usepackage{graphicx}
\usepackage{epstopdf}
\usepackage{algorithmic}
\usepackage{amsmath}
\usepackage{amssymb}
\usepackage{hyperref}
\usepackage[T1]{fontenc}
\usepackage{soul,color}
\usepackage{lineno}
\newtheorem{theorem}{Theorem}[section]
\newtheorem{proposition}[theorem]{Proposition}

\newtheorem{remark}[theorem]{Remark}
\theoremstyle{example}

\numberwithin{equation}{section}
\begin{document}


\begin{frontmatter}

%

\title{A Skew-Symmetric Energy Stable Almost Dissipation Free Formulation of the Compressible Navier-Stokes Equations}

\author[sweden,southafrica]{Jan Nordstr\"{o}m}
\cortext[secondcorrespondingauthor]{Corresponding author}
\ead{jan.nordstrom@liu.se}
\address[sweden]{Department of Mathematics, Applied Mathematics, Link\"{o}ping University, SE-581 83 Link\"{o}ping, Sweden}
\address[southafrica]{Department of Mathematics and Applied Mathematics, University of Johannesburg, P.O. Box 524, Auckland Park 2006, Johannesburg, South Africa}

\begin{abstract}
We show that a specific skew-symmetric formulation of the nonlinear terms in the compressible Navier-Stokes equations leads to an energy rate in terms of surface integrals only. No dissipative volume integrals  contribute to the energy rate.   We also discuss boundary conditions that bounds the surface integrals.
\end{abstract}

\begin{keyword}
Initial boundary value problems \sep skew-symmetric formulation \sep compressible Navier-Stokes equations \sep  energy stability  \sep summation-by-parts  \sep nonlinear boundary conditions
\end{keyword}


\end{frontmatter}


\section{Introduction}

Stability of  numerical approximations to nonlinear initial boundary value problems (IBVPs) continuous to be a challenge. The energy method \cite{kreiss1970,gustafsson1995time} applied to linear IBVPs provides energy estimates in an $L_2 $ equivalent norm together with well posed boundary conditions \cite{nordstrom2020,nordstrom2005}.
 These estimates normally require certain symmetry properties of the matrices involved. For nonlinear IBVPs, these requirements limit the use of the energy method and promted interest in entropy and kinetic energy estimates. Entropy estimates \cite{harten1983,Tadmor2003} occasionally  provide $L_2 $ estimates and boundary conditions  \cite{nordstrom2022linear}, kinetic energy estimates \cite{Jameson2008188,PIROZZOLI20107180} do not.



We have in a series of papers \cite{nordstrom2022linear-nonlinear,Nordstrom2022_Skew_Euler,NORDSTROM2024_BC} addressed the difficulties that limit the use of the energy method for nonlinear IBVPs. One difficulty relates to the nonlinear terms which must be on skew-symmetric form for Greens formula to be applicable. To address this issue, the shallow water and compressible Euler equations were transformed to skew-symmetric form in \cite{nordstrom2022linear-nonlinear,Nordstrom2022_Skew_Euler}. 
Another difficulty relates to the nonlinear boundary terms which was considered in \cite{NORDSTROM2024_BC} where energy stable boundary conditions were derived.

Once the continuous analysis leading to an energy bound is completed, one can mimic the continuous analysis discretely \cite{nordstrom_roadmap} if one can formulate the numerical procedure on summation-by-parts form
 \cite{svard2014review,fernandez2014review}.
In \cite{nordstrom2022linear-nonlinear,Nordstrom2022_Skew_Euler,NORDSTROM2024_BC}, the influence of the skew-symmetric transformation on the viscous parabolic terms were disregarded.
 In this short note we include these terms and analyse the compressible Navier-Stokes equations. 
\section{The skew-symmetric form of the compressible Euler equations}\label{sec:skew-eq}
The two-dimensional compressible Euler equations in primitive form using the ideal gas law is 
\begin{equation}\label{eq:primitive}
\begin{aligned}
  \rho_t + u_i(\rho)_{x_i} +  \rho (u_i)_{x_i}    &=0\\
    u_t + u_i (u_1)_{x_i} + p_{x_1} /\rho          &=0\\
    v_t + u_i (u_2)_{x_i}  + p_{x_2}/ \rho         &=0\\
    p_t + u_i(p)_{x_i}  + \gamma p (u_i)_{x_i} &=0\\
\end{aligned}
\end{equation}
where  $p$ is the pressure, $\gamma$ is the ratio of specific heats, $\rho$ is the density, $u_1,u_2$ are the velocities in the $(x_1,x_2)$ direction respectively. Einsteins summation convention is used above and will continue to be used below. To obtain a skew-symmetric formulation \cite{Nordstrom2022_Skew_Euler}, we choose the new dependent variables 
\begin{equation}\label{new_var}
\Phi=(\phi_1, \phi_2, \phi_3, \phi_4)^T=(\sqrt{\rho}, \sqrt{\rho} u_1, \sqrt{\rho} u_2, \sqrt{p})^T.
\end{equation}
Repeated use of the chain rule and the introduction of the new variables $\Phi$ in (\ref{new_var}) leads to the equations
\begin{equation}\label{eq:skew-symmetric-matrix-vector}
\Phi_t + B_i \Phi_{x_i} =0
\end{equation}
where $B_1,B_2$ were given in  \cite{Nordstrom2022_Skew_Euler}.
To arrive at a skew-symmetric form 
we transformed (\ref{eq:skew-symmetric-matrix-vector}) by solving
\begin{equation}\label{eq:x-direction_ODE}
(\tilde A_1  \Phi)_{x_1} + \tilde A_1^T  \Phi_{x_1} = 2 P B_1  \Phi_{x_1},  \quad (\tilde A_2  \Phi)_{x_2} + \tilde A_2^T  \Phi_{x_2} = 2 P B_2  \Phi_{x_2}
\end{equation}
for matrices $\tilde A_1, \tilde A_2$ and $P$. The solution became $P=diag(\alpha^2,  \frac{(\gamma -1)}{2},  \frac{(\gamma -1)}{2}, 1)$ with $\alpha^2>0$ arbitrary and
 \begin{equation}\label{eq:new_matrix_tildeA_4}
  \tilde A_1 = \begin{bmatrix}
      \alpha^2 u_1                                    &  &  0 & 0 \\
                                                           &  \frac{(\gamma -1)}{2}   u_1         &  0                   &   0                \\
       0                        &   0        &  \frac{(\gamma -1)}{2}  u_1                     &  0             \\
                            &  2 (\gamma -1)  \frac{\phi_4}{\phi_1}            &     0                   & (2-\gamma)u_1     
                   \end{bmatrix}, \,\
\tilde A_2 = \begin{bmatrix}
      \alpha^2 u_2                                    &  0&  0 & 0\\
      0                          &  \frac{(\gamma -1)}{2}   u_2         &   0                     &   0                \\
     0                        &  0          &      \frac{(\gamma -1)}{2}   u_2               &  0             \\
     0                    &  0       &    2 (\gamma -1)  \frac{\phi_4}{\phi_1}                    & (2-\gamma)u_2         
                   \end{bmatrix}.
 \end{equation}
The matrices $\tilde A_1, \tilde A_2$ depend on the norm matrix $P$. 
The final governing equations can be transformed to
\begin{equation}\label{eq:stab_eq}
\Phi_t + (C_1  \Phi)_{x_1} + C_2  \Phi_{x_1}+ (D_1  \Phi)_{x_2} + D_2  \Phi_{x_2}=0,
\end{equation}
where we have multiplied (\ref{eq:x-direction_ODE}) with $(2P)^{-1}$  which remove the dependence of the norm $P$ (and $\alpha^2$)  \cite{Nordstrom2022_Skew_Euler}.
\section{The skew-symmetric form of the compressible Navier-Stokes equations}\label{viscous-terms}
By adding the viscous fluxes to (\ref{eq:primitive}) we obtain
\begin{equation}\label{eq:primitive_NS}
\begin{aligned}
    \rho_t + u_i(\rho)_{x_i} +  \rho (u_i)_{x_i} &=0\\
    u_t + u_i (u_1)_{x_i} + p_{x_1} /\rho                        &=S_2= (\tau_{1j})_{x_j}/\rho\\
    v_t + u_i (u_2)_{x_i}  + p_{x_2}/ \rho                         &=S_3=(\tau_{2j})_{x_j}/\rho\\
    p_t + u_i(p)_{x_i}  + \gamma p (u_i)_{x_i}    &=S_4=(\gamma-1)( (\kappa T_{x_j})_{x_j}+ \Psi)\\
\end{aligned}
\end{equation}
where $\tau_{ij}=(\mu((u_i)_{x_j}+(u_j)_{x_i})+\lambda \delta_{ij} (u_k)_{x_k})$ is the stress tensor involving the first  viscosity $\mu$ , the second viscosity $\lambda$ and the Kronecker delta  $\delta_{ij}$ respectively. Furthermore,  $T$ is the temperature, $\kappa$ is the thermal conductivity and $\Psi= \tau_{ij}(u_i)_{x_j}$ is the positive semi-definite dissipation function \cite{White}.

To obtain a skew-symmetric version of the Navier-Stokes formulation (\ref{eq:primitive_NS}) we insert the new dependent variables in (\ref{new_var}) and obtain a new version of (\ref{eq:primitive_NS}) with coefficients multiplying $S_2,S_3,S_4$. The new set of equations have the lefthand side of (\ref{eq:skew-symmetric-matrix-vector}) and a scaled righthand side of (\ref{eq:primitive_NS}). Explicitly we get
\begin{equation}\label{eq:skew-symmetric-matrix-vector_NS}
\Phi_t + B_i \Phi_{x_i}=\Lambda S,  \quad S=(0,S_2,S_3,S_4)^T, \quad \Lambda =diag(1/(2\phi_1),\phi_1,\phi_1, 1/(2\phi_4)).
\end{equation}
Next we use the  transformation of the matrices to skew-symmetric form in (\ref{eq:x-direction_ODE}) leading to 
\begin{equation}\label{eq:stab_eq_NS_energy}
2 P \Phi_t + (\tilde A_i  \Phi)_{x_i} + \tilde A_i^T  \Phi_{x_i} =2 P \Lambda S.
\end{equation}
The new rescaled righthand side of $2 P \Lambda S$ in (\ref{eq:stab_eq_NS_energy}) is explicitly
\begin{equation}\label{RHS_NS}
2 P \Lambda S=
\begin{bmatrix}
      \alpha^2/\phi_1  &              0                        &  0 & 0 \\
                        0         &  (\gamma -1) \phi_1      &  0                   &   0                \\
       0                        &   0        &   (\gamma -1) \phi_1                      &  0             \\
                 0           & 0 &     0                   & 1/ \phi_4   
\end{bmatrix}
\begin{bmatrix}
      0      \\
      S_2  \\
       S_3  \\
       S_4
\end{bmatrix}
=(\gamma -1)
\begin{bmatrix}
    0\\
    (\tau_{1j})_{x_j}/ \phi_1 \\
    (\tau_{2j})_{x_j}/ \phi_1 \\
    ((\kappa T_{x_j})_{x_j}+ \Psi)/ \phi_4
\end{bmatrix}.
\end{equation}
\begin{remark}
The formulation (\ref{eq:stab_eq_NS_energy}) can also be obtained by starting with the conservative Navier-Stokes equations and transforming them into the primitive version which is the starting point in \cite{Nordstrom2022_Skew_Euler} and this note. All transformations up to the final skew-symmetric formulation are exact. No approximations are involved..
\end{remark}
We apply the energy method by multiplying (\ref{eq:stab_eq_NS_energy}) with  $\Phi^T$  from the left and integrate over the domain $\Omega$. Greens formula on the lefthand side and the relations  $u_1=\phi_2/ \phi_1, u_2=\phi_3/ \phi_1$ on the righthand side lead to
\begin{equation}\label{eq:boundaryPart_NS_0}
\frac{d}{dt}\|\Phi\|^2_P + \oint\limits_{\partial\Omega}  \Phi^T (n_j \tilde A_j) \Phi ds= (\gamma -1)\int\limits_{\Omega}   u_i (\tau_{ij})_{x_j} +(\kappa T_{x_j})_{x_j}  +  \Psi d{\Omega}
\end{equation}
where $(n_1,n_2)^T$ is the outward pointing unit normal from the boundary $\partial\Omega$. 
\begin{remark}
The viscous terms in (\ref{eq:boundaryPart_NS_0}) obtain a suitable form for Greens formula by the scaling with $\phi_1$.
\end{remark}
Now we can combine the two terms including the stress tensor to one complete derivative as  $u_i (\tau_{ij})_{x_j} +\Psi= u_i (\tau_{ij})_{x_j} + \tau_{ij}(u_i)_{x_j}=(u_i \tau_{ij})_{x_j}$.
By once more using Greens formula on the two remaining volume terms in (\ref{eq:boundaryPart_NS_0}) yields the final energy rate
\begin{equation}\label{eq:boundaryPart_NS_final}
\frac{d}{dt}\|\Phi\|^2_P + \oint\limits_{\partial\Omega}  \Phi^T (n_j \tilde A_j) \Phi -  (\gamma -1)(u_i \tau_{ij} n_j +\kappa T_{x_j} n_j) ds=0.
\end{equation}
\begin{remark}
At first glance, the result (\ref{eq:boundaryPart_NS_final}) is surprising. No dissipation in the form of volume integrals remains. However, a second closer look at the exact details of the norm $\|\Phi\|^2_P$ reveals that it can be seen as a linear combination of the mass and energy equations in the conservative Navier-Stokes equations.
\end{remark}

\subsection{The boundary conditions}\label{BCs}
Following \cite{NORDSTROM2024_BC}, we rotate the inviscid boundary terms into the boundary normal and tangential directions by introducing the normal $u_n=n_1u_1+n_2 u_2$ and tangential $u_t =-n_2u_1+n_1u_2$ velocities.
The rotated solution (indicated by superscript r) is $\Phi^r=(\phi_1, \phi_2^r, \phi_3r,  \phi_4)^T=(\phi_1, \phi_1 u_n, \phi_1 u_t,  \phi_4)^T$.  

Next, we  relate the temperature gradient to the normal and tangential boundary directions. We get
 \begin{equation}\label{Temp_grad}
(\gamma -1)\kappa T_{x_j} n_j= (\gamma -1)\kappa \dfrac{\partial T }{\partial n}= (\gamma -1)\kappa \dfrac{\partial  }{\partial n} \left(\dfrac{p}{R \rho}\right)=\dfrac{ (\gamma -1)\kappa  }{R} \dfrac{\partial  }{\partial n} \left(\dfrac{\phi_4}{\phi_1}\right)^2= 
2 \theta  \left(\dfrac{\phi_4}{\phi_1}\right) \psi_n
\end{equation}
where  $\theta =  (\gamma -1)\kappa/R=\gamma \mu/Pr$,  $Pr$ is the Prandtl number and $\psi_n$ is short notation for $\partial (\phi_4/\phi_1)/\partial n$ . The stress tensor related term $u_i \tau_{ij} n_j $ given in the  $(x_1,x_2)$ coordinates is expanded as
 \begin{equation}\label{Stress_1}
 u_i \tau_{ij} n_j = u_1 (\tau_{11} n_1 + \tau_{12} n_2)+ u_2  (\tau_{21} n_1 + \tau_{22} n_2)=u_1 \tau_1+u_2  \tau_2= \vec u  \cdot \vec \tau.
\end{equation}
The velocity $ \vec u$ and stress vector  $\vec \tau$ in the normal and tangential directions are expressed as
 \begin{equation}\label{Stress_2}
 \vec u_r=
 \begin{bmatrix}
      u_n  \\
      u_t                            
\end{bmatrix}=N  \vec u 
, \quad 
 \vec \tau_r=
 \begin{bmatrix}
       \tau_n  \\
       \tau_t                            
\end{bmatrix}=N  \vec \tau
 \quad   \text{where} \quad  
N=
 \begin{bmatrix}
       n_1  & n_2\\
       -n_2    & n_1                     
\end{bmatrix}
\quad   \text{and} \quad  N^T N= I_2=diag(1,1).
\end{equation}
Using (\ref{Stress_1}) and  (\ref{Stress_2}) we obtain the rotated stress tensor related term as  
 \begin{equation}\label{Stress_3}
 u_i \tau_{ij} n_j =\vec u  \cdot \vec \tau=(N^T \vec u_r)  \cdot (N^T \vec \tau_r)=\vec u_r^T N N^T \vec \tau_r=\vec u_r \cdot \vec \tau_r.
 \end{equation}
\begin{remark}
The viscous terms add three new dependent variables ($\tau_n,\tau_t$ and $\psi_n$) that involve gradients.
\end{remark}

We can now combine the seven variables ($\Phi^r$ and $\psi_n,\tau_n,\tau_t$) in the surface integral of (\ref{eq:boundaryPart_NS_final}) into a vector-matrix-vector product $BT=W^T A W/{w_1}$, where $W=(w_1,w_2,w_3)^T$, $w_1=\phi_1, w_2=(\phi_2^r, \phi_3^r,  \phi_4)^T$ and $w_3=(\theta\psi_n, \frac{(\gamma -1)}{2} \tau_n, \frac{(\gamma -1)}{2} \tau_t)^T$. The block matrix  $A$  in $BT$  is given by 
\begin{equation}\label{boundarmatrix_contraction_rotated_expanded}
A=
\begin{bmatrix}
    A_{11}&  0         &   0    \\
      0     & A_{22} &  -I_3 \\
      0    &-I_3      &   0       
 \end{bmatrix},   \quad  A_{11}  = \phi_2^r,  \quad 
   A_{22}=
 \begin{bmatrix}
    \frac{(\gamma -1)}{2} \phi_2^r  & 0 &   (\gamma -1)   \phi_4   \\
      0     &  \frac{(\gamma -1)}{2} \phi_2^r&  0\\
       (\gamma -1)   \phi_4       & 0     &   (2-\gamma) \phi_2^r         
 \end{bmatrix}.
\end{equation}
The zero matrices in (\ref{boundarmatrix_contraction_rotated_expanded}) have appropriate sizes and $ I_3 = diag(1,1,1)$. 

Following \cite{NORDSTROM2024_BC} we now transform  the boundary term in (\ref{boundarmatrix_contraction_rotated_expanded}) into diagonal form.  By introducing the transformation matrix $R$ we obtain $BT$ with block diagonal matrix $\Lambda$ as
\begin{equation}\label{blockrotate}
BT=
\frac{1}{w_1}  (RW)^T \Lambda_A (RW)= \frac{1}{w_1}
\begin{bmatrix}
  A_{11} w_1\\
  A_{22} w_2-w_3\\
  w_3
\end{bmatrix}^T
\begin{bmatrix}
    A_{11}^{-1}&  0         &   0    \\
      0     & A_{22}^{-1}&  0 \\
      0    &    0  & -A_{22}^{-1} 
 \end{bmatrix}
 \begin{bmatrix}
  A_{11} w_1\\
  A_{22} w_2-w_3\\
  w_3
\end{bmatrix}.
 \end{equation}
Since $ A_{11}$ is a scalar we continue by diagonalising $A_{22}$ using the transformation matrix $S_{22}$ which lead to 
\begin{equation}\label{diagrotate}
A_{22}^{-1}=
S_{22} \Lambda_{22}^{-1} S_{22}^T=
\begin{bmatrix}
    1&  0         &   -2  \frac{\phi_4}{\phi_1}  \\
      0     & 1&  0 \\
      0    &    0  & 1
 \end{bmatrix}
\begin{bmatrix}
      \frac{(\gamma -1)}{2} \phi_2^r  &  0         &   0    \\
      0     & \frac{(\gamma -1)}{2} \phi_2^r &  0 \\
      0    &    0  & (2-\gamma) \phi_2^r  \beta(M_n^2)
 \end{bmatrix}^{-1}
 \begin{bmatrix}
    1&  0         &  0   \\
      0     & 1&  0 \\
       -2  \frac{\phi_4}{\phi_1}     &    0  & 1
 \end{bmatrix}.
 \end{equation}
In (\ref{diagrotate}) we have  $\beta=(M_n^2- \frac{2(\gamma -1)}{\gamma(2-\gamma)})/M_n^{2}$,  $M_n=u_n/c$ is the normal Mach number and $c$ is the speed of sound. By introducing $S=diag(1,S_{22},S_{22})$ 
 we find the final result
 \begin{equation}\label{totalrotate_3}
 BT= (S^T R W)^T  \Lambda (S^T R W),
\quad  \text{where}  \quad 
 \Lambda
=\frac{1}{\phi_1^2 u_n}
\begin{bmatrix}
  1 &  0         &   0 \\
      0     & \mathsf{\Lambda}  &  0 \\
      0    &    0  & -\mathsf{\Lambda} 
 \end{bmatrix} , 
\quad 
\mathsf{\Lambda} =
\begin{bmatrix}
    \frac{2}{\gamma -1}   &  0         &   0    \\
      0     &    \frac{2}{\gamma -1}  \\
      0    &     0 &    \frac{1}{(2-\gamma )\beta(M_n^2)}  
 \end{bmatrix}.
 \end{equation}
The number of boundary conditions is given by the number of negative entries in $\Lambda$ and varies with the sign of $u_n$ and the size of $M_n^2$.
Complete details for strong and weak boundary conditions are given in \cite{NORDSTROM2024_BC}.
 \begin{remark}
With the boundary terms limited by appropriate boundary conditions, we will obtain an estimate in terms of $\|\Phi\|^2_P$ which is a weaker norm (not including gradients) than expected from an incompletely parabolic initial boundary value problem.
\end{remark}


\section{Summary and outlook}\label{sec:conclusion}
We have shown that a specific skew-symmetric formulation of the compressible Navier-Stokes equations leads to an energy rate in terms of surface integrals only. 
We also briefly discussed boundary conditions. With bounded surface integrals, the resulting energy estimate is obtained in a weaker norm (not including gradients) than expected from an incompletely parabolic initial boundary value problem. In future work we will derive explicit boundary conditions and prove nonlinear semi-discrete energy stability.

\section*{Acknowledgments}

J. N. was supported by Vetenskapsr{\aa}det, Sweden [award 2021-05484 VR] and University of Johannesburg.

\bibliographystyle{elsarticle-num}
\bibliography{References_Jan}

\begin{thebibliography}{10}
\expandafter\ifx\csname url\endcsname\relax
  \def\url#1{\texttt{#1}}\fi
\expandafter\ifx\csname urlprefix\endcsname\relax\def\urlprefix{URL }\fi
\expandafter\ifx\csname href\endcsname\relax
  \def\href#1#2{#2} \def\path#1{#1}\fi

\bibitem{kreiss1970}
H.-O. Kreiss, Initial boundary value problems for hyperbolic systems, Commun.
  Pur. Appl. Math. 23~(3) (1970) 277--298.

\bibitem{gustafsson1995time}
B.~Gustafsson, H.-O. Kreiss, J.~Oliger, Time dependent problems and difference
  methods, Vol.~24, JWS, 1995.

\bibitem{nordstrom2020}
J.~Nordstr\"{o}m, T.~M. Hagstrom, The number of boundary conditions for initial
  boundary value problems, SIAM Journal on Numerical Analysis 58~(5) (2020)
  2818--2828.

\bibitem{nordstrom2005}
J.~Nordstr\"{o}m, M.~Sv\"{a}rd, Well posed boundary conditions for the
  {N}avier--{S}tokes equations, SIAM J. Numer. Anal. 43 (2005) 1231--1255.

\bibitem{harten1983}
A.~Harten, On the symmetric form of systems of conservation laws with entropy,
  J. Comput. Phys. 49 (1983) 151--164.

\bibitem{Tadmor2003}
E.~Tadmor, Entropy stability theory for difference approximations of nonlinear
  conservation laws and related time-dependent problems, Acta Numer. 12 (2003)
  451--512.

\bibitem{nordstrom2022linear}
J.~Nordström, A.~R. Winters, A linear and nonlinear analysis of the shallow
  water equations and its impact on boundary conditions, Journal of
  Computational Physics 463 111254 (2022).

\bibitem{Jameson2008188}
A.~Jameson, Formulation of kinetic energy preserving conservative schemes for
  gas dynamics and direct numerical simulation of one-dimensional viscous
  compressible flow in a shock tube using entropy and kinetic energy preserving
  schemes, Journal of Scientific Computing 34~(2) (2008) 188 – 208.

\bibitem{PIROZZOLI20107180}
S.~Pirozzoli, Generalized conservative approximations of split convective
  derivative operators, Journal of Computational Physics 229~(19) (2010)
  7180--7190.

\bibitem{nordstrom2022linear-nonlinear}
J.~Nordström, Nonlinear and linearised primal and dual initial boundary value
  problems: When are they bounded? how are they connected?, Journal of
  Computational Physics 455 111001 (2022).

\bibitem{Nordstrom2022_Skew_Euler}
J.~Nordström, A skew-symmetric energy and entropy stable formulation of the
  compressible {E}uler equations, Journal of Computational Physics 470 111573
  (2022).

\bibitem{NORDSTROM2024_BC}
J.~Nordström, Nonlinear boundary conditions for initial boundary value
  problems with applications in computational fluid dynamics, Journal of
  Computational Physics 498 (2024) 112685.

\bibitem{nordstrom_roadmap}
J.~Nordstr\"{o}m, A roadmap to well posed and stable problems in computational
  physics, {J. Sci. Comput.} 71~(1) (2017) 365--385.

\bibitem{svard2014review}
M.~Sv{\"a}rd, J.~Nordstr{\"o}m, Review of summation-by-parts schemes for
  initial--boundary-value problems, Journal of Computational Physics 268 (2014)
  17--38.

\bibitem{fernandez2014review}
D.~C. Del Rey~Fern{\'a}ndez, J.~E. Hicken, D.~W. Zingg, Review of
  summation-by-parts operators with simultaneous approximation terms for the
  numerical solution of partial differential equations, Computers \& Fluids 95
  (2014) 171--196.

\bibitem{White}
F.~M. White, Viscous Fluid Flow 2nd ed., McGraw-Hill Inc., 1974.

\end{thebibliography}

\end{document}